\documentclass{svjour3}                     
\smartqed  
\usepackage{graphicx}
\usepackage{amssymb,amsmath,amsfonts,latexsym,euscript}
\def\d{\,{\rm d}}

\begin{document}

\title{Construction of interpolation splines minimizing the semi-norm in the space $K_2(P_m)$}


\titlerunning{Construction of interpolation splines minimizing the semi-norm in the space $K_2(P_m)$ }        

\author{Abdullo R.~Hayotov}


\institute{
Abdullo Rakhmonovich Hayotov \at
Institute of Mathematics, National University of Uzbekistan, Do`rmon yo`li str.,29, Tashkent, 100125, Uzbekistan\\
              \email{hayotov@mail.ru}
}

\date{Received: date / Accepted: date}

\maketitle

\begin{abstract}

In the present paper, using S.L. Sobolev's method, interpolation
 spline that minimizes the expression
$\int_0^1(\varphi^{(m)}(x)+\omega^2\varphi^{(m-2)}(x))^2dx$ in the
 $K_2(P_m)$ space
are constructed. Explicit formulas for the coefficients of the
interpolation splines are obtained. The obtained interpolation
spline is exact for monomials $1,x,x^2,..., x^{m-3}$ and for trigonometric functions $\sin\omega x$ and $\cos\omega x$.

\keywords{Interpolation spline \and Hilbert space \and the norm
minimizing property \and S.L.~Sobolev's method \and discrete
argument function}
 \subclass{MSC 41A05, 41A15 }
\end{abstract}

\section{Introduction. Statement of the Problem}

In order to find an approximate representation of a function
$\varphi$ by elements of a certain finite dimensional space, it is
possible to use values of this function at some finite set of
points $x_\beta$,\ $\beta=0,1,...,N$. The corresponding problem is
called \emph{the interpolation problem}, and the points $x_\beta$
are called \emph{the interpolation nodes}.

There are polynomial and spline interpolations. It is known that
the polynomial approximation is non-practical for approximation of
functions with finite and small smoothness, which often occurs in
applications. This circumstance makes necessary to work with the
splines. Spline functions are very useful in applications. Classes
of spline functions possess many nice structural properties as
well as excellent approximation powers. They are used, for
example, in data fitting, function approximation, numerical
quadrature, and the numerical solution of ordinary and partial
differential equations, integral equations, and so on. Many books
are devoted to the theory of splines, for example, Ahlberg et al
\cite{Ahlb67}, Arcangeli et al \cite{Arc04}, Attea \cite{Attea92},
Berlinet and Thomas-Agnan \cite{BerAgnan04}, Bojanov et al
\cite{Boj93}, de Boor \cite{Boor78}, Eubank \cite{Eub88}, Green
and Silverman \cite{GrSi94}, Ignatov and Pevniy \cite{Ign91},
Korneichuk et al \cite{Korn93}, Laurent \cite{Lor75}, Mastroianni
and Milovanovi\'c \cite{GM_GVM}, N\"{u}rnberger \cite{Nur89},
Schumaker \cite{Schum81},  Stechkin and Subbotin \cite{Stech76},
Vasilenko \cite{Vas83}, Wahba \cite{Wahba90} and others.

If the exact values $\varphi(x_\beta)$ of an unknown smooth
function $\varphi(x)$ at the set of points $\{x_\beta,\
\beta=0,1,...,N\}$ in an interval $[a,b]$ are known, it is usual
to approximate $\varphi$ by minimizing
\begin{equation}
\label{eq.(1.1)} \int_a^b(g^{(m)}(x))^2dx
\end{equation}
in the set of interpolating functions (i.e.,
$g(x_\beta)=\varphi(x_\beta)$, $\beta=0,1,...,N$) of the space $L_2^{(m)}(a,b)$. Here $L_2^{(m)}(a,b)$ is the Sobolev space
of functions with a square integrable $m$-th generalized
derivative. It turns out that the solution is a natural polynomial
spline of degree $2m-1$ with knots $x_0,x_1,...,x_N$ called
\emph{the interpolating $D^m$-spline} for the points
$(x_\beta,\varphi(x_\beta))$. In the non periodic case this
problem has been investigated, at the first time, by Holladay
\cite{Hol57} for $m=2$. His results have been generalized by de
Boor \cite{deBoor63} for any $m$. In the Sobolev space
$\widetilde{L_2^{(m)} }$ of periodic functions, the minimization
problem of integrals of type (\ref{eq.(1.1)}) was investigated in
works \cite{Freeden84,Freeden88,Golomb68,MamHayShad13,Schoen64}
and others.

We consider the Hilbert space
\[
K_2(P_m)=\Bigl\{\varphi:[0,1]\to \mathbb{R}\ \Bigm|\
\varphi^{(m-1)} \mbox{ is absolutely continuous and }
\varphi^{(m)}\in L_2(0,1)\Bigr\},
\]
equipped with the norm
\begin{equation}\label{eq.(1.2)}
\left\| {\varphi \,|K_2(P_m)} \right\| = \left\{ {\int_0^1 {\left(
{P_m \left( {{\d \over {\d x}}} \right)\varphi (x)} \right)} ^2 \d
x} \right\}^{1/2},
\end{equation}
where
\[
P_m \left( {{\d \over {\d x}}} \right)=\frac{\d^m}{\d
x^m}+\omega^2\frac{\d^{m-2}}{\d x^{m-2}},\ \
\omega>0,\ \ \  m\geq 2
\]
and
\[
{\int_0^1 {\left( {P_m \left( {{\d
\over {\d x}}} \right)\varphi (x)} \right)} ^2 \d x}<\infty.\] The
equality (\ref{eq.(1.2)}) is the semi-norm and $\|\varphi\|=0$ if
and only if $\varphi(x)=c_1\sin \omega x+c_2\cos\omega x+R_{m-3}(x)$, where
$R_{m-3}(x)$ is a polynomial of degree $m-3$.

It should be noted that for a linear differential operator of
order $n$, $L\equiv P_n(\d/\d x)$, Ahlberg,  Nilson,  and Walsh in
the book \cite[Chapter 6]{Ahlb67} investigated the Hilbert spaces
in the context of generalized splines. Namely, with the inner
product
\[\langle\varphi,\psi\rangle=\int_0^1 L\varphi (x)\cdot L\psi (x)  \d x ,\]
$K_2(P_n)$ is a Hilbert space if we identify functions that differ by a solution of $L\varphi=0$.

Consider the following interpolation problem:

\begin{problem}\label{Pr1}
To find the function $S_m(x)\in K_2(P_m)$, which gives the
minimum of the norm (\ref{eq.(1.2)}) and satisfies the
interpolation condition
\begin{equation}\label{eq.(1.3)}
S_m(x_\beta)=\varphi(x_\beta),\ \ \beta=0,1,...,N,
\end{equation}
where $x_\beta\in [0,1]$  are the nodes of interpolation, $\varphi(x_\beta)$ are given values.
\end{problem}

Following \cite[p.46, Theorem 2.2]{Vas83}  we get the analytic
representation of the interpolation spline $S_m(x)$
\begin{equation}\label{eq.(1.4)}
S_m(x)=\sum\limits_{\gamma=0}^NC_\gamma
G_m(x-x_\gamma)+d_1\sin(\omega x)+d_2\cos(\omega x)+R_{m-3}(x),
\end{equation}
where $C_\gamma$, $\gamma=0,1,...,N$, $d_1$ and $d_2$  are real numbers,
$R_{m-3}(x)=\sum\limits_{\alpha=0}^{m-3}r_{\alpha}x^{\alpha}$  is a
polynomial of degree $m-3$ and
\begin{equation}\label{eq.(1.5)}
G_m(x) = \frac{(-1)^m \mathrm{sign}x}{4\omega^{2m-1}}\,\left(
(2m-3)\sin \omega x-\omega x \cos \omega x
+2\sum\limits_{k=1}^{m-2}\frac{(-1)^k(m-k-1)(\omega
x)^{2k-1}}{(2k-1)!}\right)
\end{equation}
is a fundamental solution of the operator
$\frac{\d^{2m}}{\d x^{2m}}+\omega^2 \frac{\d^{2m-2}}{\d x^{2m-2}}+\omega^4\frac{\d^{2m-4}}{\d x^{2m-4}}$, i.e., $G_m(x)$ is a solution of the
equation
\begin{equation}\label{eq.(1.6)}
G_m^{(2m)} (x) + 2\omega^2G_m^{(2m-2)}(x) +\omega^4
G_m^{(2m-4)}(x) = \delta(x),
\end{equation}
here $\delta(x)$ is Dirac's delta function.

It is known that (see, for instance,
\cite{Vas83}) the solution $S_m(x)$ of
the form (\ref{eq.(1.4)}) of Problem \ref{Pr1}  exists, is unique
when $N+1\geq m$ and coefficients $C_\gamma$, $d_1$, $d_2$ and $r_\alpha$ of
$S_m(x)$ are defined by the following system of $N+m+1$ linear
equations
\begin{eqnarray}
&&\sum\limits_{\gamma=0}^NC_{\gamma}G_m(x_\beta-x_\gamma)+d_1\sin(\omega x_\beta)+d_2\cos(\omega x_\beta)+R_{m-3}(x_\beta)=\varphi(x_\beta),\ \ \label{eq.(1.7)}\\
&&\beta=0,1,...,N,\nonumber \\
&& \sum\limits_{\gamma=0}^NC_{\gamma}\sin(\omega x_\gamma)=0,\label{eq.(1.8)}\\
&& \sum\limits_{\gamma=0}^NC_{\gamma}\cos(\omega x_\gamma)=0,\label{eq.(1.9)}\\
&& \sum\limits_{\gamma=0}^NC_{\gamma}x_\gamma^\alpha=0,\ \
\alpha=0,1,...,m-3.\label{eq.(1.10)}
\end{eqnarray}

The main aim of the present paper is to solve Problem 1, i.e., to
solve system (\ref{eq.(1.7)})-(\ref{eq.(1.10)}) for equally
spaced nodes $x_\beta=h\beta$, $\beta=0,1,...,N,$ $h=1/N$,
$N+1\geq m$ and to find analytic formulas for the coefficients
$C_\gamma$, $d_1$, $d_2$ and $r_\alpha$ of $S_m(x)$.

It should be noted that, using Sobolev method, interpolation
splines minimizing the semi-norms in the $L_2^{(m)}(0,1)$, $W_2^{(m,m-1)}(0,1)$ and $K_2(P_2)$ Hilbert spaces were constructed in works
\cite{CabHaySha14,Hay14,HayMilShad13,HayMilShad14,ShadHay12,ShadHayAzam14}.
Furthermore connection between interpolation spline  and optimal
quadrature formula in the sense of Sard in $L_2^{(m)}(0,1)$ and $K_2(P_2)$ spaces were
shown in \cite{CabHaySha14} and \cite{HayMilShad13}.

The rest of the paper is organized as follows: in Section 2 we
give  some definitions and known results. In Section 3 it is given
the algorithm for solution of system
(\ref{eq.(1.7)})-(\ref{eq.(1.10)}) when the nodes $x_\beta$ are
equally spaced. Using this algorithm, the coefficients of the
interpolation spline $S_m(x)$ are computed in Section 4.

\section{Preliminaries}
\setcounter{equation}{0} \setcounter{theorem}{0}
\setcounter{definition}{0} \setcounter{lemma}{0}
\setcounter{example}{0} \setcounter{figure}{0}
\setcounter{table}{0}

In this section we give some definitions and known results that we
need to prove the main results.

Below mainly we use the concept of discrete argument functions and
operations on them. The theory of discrete argument functions is
given in \cite{Sob74,SobVas}. For completeness we give some
definitions about functions of discrete argument.

Assume that the nodes $x_\beta$ are equal spaced, i.e., $x_\beta=
h\beta,$ $h = {1 \over N}$, $N = 1,2,...$.

\textbf{Definition 2.1.}{ The function $\varphi (h\beta )$ is a}
{\it function of discrete argument}  {if it is given on some set
of integer values of} $\beta$.

\textbf{Definition 2.2.} {\it The inner product} {of two discrete
functions $\varphi(h\beta )$ and $\psi (h\beta )$ is given by}
$$
\left[ {\varphi(h\beta),\psi(h\beta) } \right] =
\sum\limits_{\beta  =  - \infty }^\infty  {\varphi (h\beta ) \cdot
\psi (h\beta )},
$$
{if the series on the right hand side of the last equality
converges absolutely.}

\textbf{Definition 2.3.} {\it The convolution} \emph{of two
functions $\varphi(h\beta )$ and $\psi (h\beta )$ is the inner
product}
$$
\varphi (h\beta )*\psi (h\beta ) = \left[ {\varphi (h\gamma ),\psi
(h\beta  - h\gamma )} \right] = \sum\limits_{\gamma  =  - \infty
}^\infty  {\varphi (h\gamma ) \cdot \psi (h\beta  - h\gamma )}.
$$

\emph{The Euler-Frobenius polynomials} $E_k (x)$, $k = 1,2,...$ is
defined by the following formula \cite{SobVas}
\begin{equation}\label{eq.(2.1)}
 E_k (x) = \frac{{(1 - x)^{k + 2} }}{x}\left( {x\frac{d}{{dx}}}
\right)^k \frac{x}{{(1 - x)^2 }},
\end{equation}
$E_0 (x) = 1$.

For the Euler-Frobenius polynomials $E_k(x)$ the following
identity holds
\begin{equation}\label{eq.(2.2)}
E_k (x) = x^k E_k \left( {\frac{1}{x}} \right),
\end{equation}
and also the following theorem is true

\begin{theorem}\label{THM:2.1}
\emph{ (Lemma 3 of \cite{Shad10})}. {Polynomial $Q_k (x)$ which is
defined by the formula}
\begin{equation}\label{eq.(2.3)}
Q_k (x) = (x - 1)^{k + 1} \sum\limits_{i = 0}^{k + 1}
{\frac{{\Delta ^i 0^{k + 1} }}{{(x - 1)^i }}}
\end{equation}
{is the Euler-Frobenius polynomial (2.1) of degree $k$, i.e.
$Q_k(x) = E_k(x)$, where
$\Delta^i0^k=\sum_{l=1}^i(-1)^{i-l}C_i^ll^k.$}
\end{theorem}

The following formula is valid \cite{Ham73}:
\begin{equation}\label{eq.(2.4)}
\sum\limits_{\gamma  = 0}^{n - 1} {q^\gamma  \gamma ^k  =
\frac{1}{{1 - q}}\sum\limits_{i = 0}^k {\left( {\frac{q}{{1 - q}}}
\right)^i \Delta ^i 0^k  - \frac{{q^n }}{{1 - q}}\sum\limits_{i =
0}^k {\left( {\frac{q}{{1 - q}}} \right)^i \Delta ^i \gamma ^k
|_{\gamma  = n} ,} } }
\end{equation}
where $\Delta^i\gamma^k$ is the finite difference of order $i$ of
$\gamma^k$, $q$ is the ratio of a geometric progression. When
$|q|<1$ from (\ref{eq.(2.4)}) we have
\begin{equation}\label{eq.(2.5)}
\sum\limits_{\gamma  = 0}^{\infty} q^\gamma  \gamma ^k  =
\frac{1}{{1 - q}}\sum\limits_{i = 0}^k \left( {\frac{q}{{1 - q}}}
\right)^i \Delta ^i 0^k.
\end{equation}
In our computations we need the discrete analogue $D_m(h\beta)$ of
the differential operator $\frac{\d^{2m}}{\d x^{2m}}+\omega^2\frac{\d^{2m-2}}{\d x^{2m-2}}+\omega^4\frac{\d^{2m-4}}{\d x^{2m-4}}$ which satisfies the
following equality
\begin{equation}\label{eq.(2.6)}
D_m(h\beta)*G_m(h\beta) = \delta(h\beta),
\end{equation}
where $G_m(h\beta)$ is the discrete argument function corresponding to
$G_m(x)$ defined by (\ref{eq.(1.5)}), $\delta(h\beta)$ is equal to
0 when $\beta \ne 0$ and is equal to 1 when $\beta = 0$, i.e.
$\delta(h\beta)$ is the discrete delta-function. The equation
(\ref{eq.(2.6)}) is the discrete analogue of the equation
(\ref{eq.(1.6)}).

In \cite{Hay13,Hay14} the discrete analogue $D_m(h\beta)$ of the
differential operator $\frac{\d^{2m}}{\d x^{2m}}+\omega^2\frac{\d^{2m-2}}{\d x^{2m-2}}+\omega^4\frac{\d^{2m-4}}{\d x^{2m-4}}$, which satisfies  equation
(\ref{eq.(2.6)}), is constructed and the following is
proved.

\begin{theorem}\label{THM2.2}
 The discrete analogue to the differential operator
$\frac{\d^{2m}}{\d x^{2m}}+2\omega^2\frac{\d^{2m-2}}{\d
x^{2m-2}}+\omega^4\frac{\d^{2m-4}}{\d x^{2m-4}}$ satisfying
equation \emph{(\ref{eq.(2.6)})} has the form
\begin{equation}\label{eq.(2.7)}
D_m(h\beta)=p\left\{
\begin{array}{ll}
\sum\limits_{k=1}^{m-1}A_k \lambda_k^{|\beta|-1},& |\beta|\geq 2,\\
1+\sum\limits_{k=1}^{m-1}A_k,& |\beta|=1,\\
C+\sum\limits_{k=1}^{m-1}\frac{A_k}{\lambda_k},& \beta=0,
\end{array}
\right.
\end{equation}
where
\begin{eqnarray}
A_k&=&\frac{(1-\lambda_k)^{2m-4}(\lambda_k^2-2\lambda_k\cos
h\omega+1)^2p_{2m-2}^{(2m-2)}}{\lambda_k{\cal P}_{2m-2}'(\lambda_k)},\label{eq.(2.8)}\\
C&=&4-4\cos h\omega-2m-\frac{p_{2m-3}^{(2m-2)}}{p_{2m-2}^{(2m-2)}},\ \ \ p=\frac{2\omega^{2m-1}}{(-1)^mp_{2m-2}^{(2m-2)}},\label{eq.(2.9)}\\
p_{2m-2}^{(2m-2)}&=&(2m-3)\sin h\omega-h\omega\cos h\omega\nonumber \\
&&\qquad\qquad\qquad
+2\sum\limits_{k=1}^{m-2}\frac{(-1)^k(m-k-1)(h\omega)^{2k-1}}{(2k-1)!},\label{eq.(2.10)}
\end{eqnarray}
$$
{\cal P}_{2m-2}(x)=\sum_{s=0}^{2m-2}p_s^{(2m-2)}x^s
=(1-x)^{2m-4}\bigg[[(2m-3)\sin h\omega-h\omega \cos h\omega]x^2
$$
$$
\qquad\qquad\qquad +[2h\omega-(2m-3)\sin(2h\omega)]x+[(2m-3)\sin
h\omega-h\omega \cos h\omega]\bigg]
$$
\begin{equation}\label{eq.(2.11)}
+2(x^2-2x\cos
h\omega+1)^2\sum_{k=1}^{m-2}\frac{(-1)^k(m-k-1)(h\omega)^{2k-1}(1-x)^{2m-2k-4}E_{2k-2}(x)}{(2k-1)!},
\end{equation}
here $E_{2k-2}(x)$ is the Euler-Frobenius polynomial of degree
$2k-2$, $\omega>0$, $h\omega\leq 1$, $h=1/N$, $N\geq m-1$, $m\geq
2$, $p_{2m-2}^{(2m-2)},\ p_{2m-3}^{(2m-2)}$ are the coefficients
and $\lambda_k$ are the roots of the polynomial ${\cal
P}_{2m-2}(\lambda)$, $|\lambda_k|<1$.
\end{theorem}

Furthermore several properties  of the discrete argument function
$D_m(h\beta)$ were given in \cite{Hay13,Hay14}. Here we give the
following properties of the discrete argument function $D_m(h\beta)$
which we need in our computations.

\begin{theorem}\label{THM2.3}
The discrete analogue  $D_m(h\beta)$ of the differential operator
$\frac{\d^{2m}}{\d x^{2m}}+2\omega^2\frac{\d^{2m-2}}{\d
x^{2m-2}}+\omega^4\frac{\d^{2m-4}}{\d x^{2m-4}}$ satisfies the
following equalities

\emph{1)} $D_m(h\beta)*\sin(h\omega\beta)=0,$

\emph{2)} $D_m(h\beta)*\cos(h\omega\beta)=0,$

\emph{3)} $D_m(h\beta)*(h\omega\beta)\sin(h\omega\beta)=0,$

\emph{4)} $D_m(h\beta)*(h\omega\beta)\cos(h\omega\beta)=0,$

\emph{5)} $D_m(h\beta)*(h\beta)^\alpha=0,$ $\alpha=0,1,...,2m-5$.
\end{theorem}

\section{The algorithm for computation of coefficients of interpolation splines}
\setcounter{equation}{0} \setcounter{theorem}{0}
\setcounter{definition}{0}
 \setcounter{lemma}{0} \setcounter{example}{0}
 \setcounter{figure}{0} \setcounter{table}{0}

In the present section we give the algorithm for solution of
system (\ref{eq.(1.7)})-(\ref{eq.(1.10)}) when the nodes $x_\beta$
are equally spaced, i.e., $x_\beta= h\beta,$ $h = {1 \over N}$, $N
= 1,2,...$ . Here we use similar method suggested by S.L. Sobolev
\cite{Sob77,SobVas} for finding the coefficients of optimal
quadrature formulas in the Sobolev space $L_2^{(m)}(0,1)$.

Suppose that $C_\beta=0$  when $\beta < 0$ and $\beta  > N$. Using
Definition 2.3, we rewrite system
(\ref{eq.(1.7)})-(\ref{eq.(1.10)}) in the convolution form
\begin{eqnarray}
&&G_m(h\beta )*C_\beta+ d_1\sin(h\omega\beta)+d_2\cos(h\omega\beta)+R_{m-3}(h\beta)=
\varphi(h\beta),\label{eq.(3.1)}\\
&&\beta  = 0,1,...,N,\nonumber \\
&&\sum\limits_{\beta=0}^N C_\beta\cdot \sin(h\omega\beta)=0,\label{eq.(3.2)}\\
&&\sum\limits_{\beta=0}^N C_\beta\cdot \cos(h\omega\beta)=0,\label{eq.(3.3)}\\
&&\sum\limits_{\beta=0}^N C_\beta\cdot (h\beta)^{\alpha}=0,\ \ \
\alpha=0,1,...,m-3, \label{eq.(3.4)}
\end{eqnarray}
where
$R_{m-3}(h\beta)=\sum_{\alpha=0}^{m-3}r_\alpha(h\beta)^\alpha$.

Thus we have the following problem.

\begin{problem}\label{Pr2}
Find the coefficients $C_\beta$, $(\beta=0,1,...,N)$, $d_1$, $d_2$  and polynomial
$R_{m-3}(h\beta)$ of degree $m-3$ which satisfy system
(\ref{eq.(3.1)})-(\ref{eq.(3.4)}).
\end{problem}

Further we investigate Problem \ref{Pr2} which is equivalent to
Problem \ref{Pr1}. Instead of $C_\beta$ we introduce the following
functions
\begin{eqnarray}
v(h\beta ) &=& G_m(h\beta )*C_\beta,\label{eq.(3.5)}\\
u(h\beta)&=&v\left({h\beta}\right)+d_1\sin(h\omega\beta)+d_2\cos(h\omega\beta)+R_{m-3}(h\beta).\label{eq.(3.6)}
\end{eqnarray}
Now we express the coefficients
$C_\beta$ by the function $u(h\beta )$.

Taking into account (\ref{eq.(2.7)}), (\ref{eq.(3.6)}) and
Theorems \ref{THM2.2}, \ref{THM2.3}, for the coefficients we
have
\begin{equation}\label{eq.(3.7)}
C_\beta=D_m(h\beta )*u(h\beta ).
\end{equation}

Thus, if we find the function $u(h\beta)$, then the coefficients
$C_{\beta}$ will be found from equality (\ref{eq.(3.7)}).

To calculate the convolution (\ref{eq.(3.7)}) it is required to
find the representation of the function $u(h\beta )$ for all
integer values of $\beta$. From equality (\ref{eq.(3.1)}) we get
that $u(h\beta )=\varphi(h\beta )$ when $h\beta \in [0,1]$. Now we
need to find the representation of the function $u(h\beta )$ when
$\beta < 0$ and $\beta>N$.

Since $C_\beta= 0$ when $h\beta \notin [0,1]$ then
$$
C_\beta = D_m(h\beta )*u(h\beta ) = 0,\,\,\,\,\,\,\,\,\,h\beta
\notin [0,1].
$$

Now we calculate the convolution $v(h\beta ) =
G_m(h\beta)*C_\beta$ when $\beta\leq 0$ and $\beta\geq N$.

Suppose $\beta  \leq 0$ then taking into account equalities
(\ref{eq.(1.5)}), (\ref{eq.(3.2)})-(\ref{eq.(3.4)}), we have
$$
\begin{array}{rcl}
 v(h\beta ) &=& \sum\limits_{\gamma  =  - \infty }^\infty
{C_\gamma
\,G_m(h\beta  - h\gamma )  }\\
&=&\sum\limits_{\gamma=0}^NC_{\gamma}\frac{(-1)^m \mathrm{sign}(h\beta-h\gamma)}{4\omega^{2m-1}}\bigg[
(2m-3)\sin(h\omega\beta-h\omega\gamma)-(h\omega\beta-h\omega\gamma)\cos(h\omega\beta-h\omega\gamma)\\
&&+2\sum\limits_{k=1}^{m-2}\frac{(-1)^k(m-k-1)(h\omega\beta-h\omega\gamma)^{2k-1}}{(2k-1)!}\Bigg]\\
&=&-\frac{(-1)^m}{4\omega^{2m-1}}\Bigg[\cos(h\omega\beta)\sum\limits_{\gamma=0}^NC_{\gamma}(h\omega\gamma)\cos(h\omega\gamma)+\sin(h\omega\beta)\sum\limits_{\gamma=0}^N
C_{\gamma}(h\omega\gamma)\sin(h\omega\gamma)\\
&&+2\sum\limits_{k=[\frac{m}{2}]}^{m-2}\sum\limits_{\alpha=m-2}^{2k-1}\frac{(-1)^{k+\alpha}(m-k-1)\omega^{2k-1}}{(2k-1-\alpha)!\alpha!}
(h\beta)^{2k-1-\alpha} \sum\limits_{\gamma=0}^NC_{\gamma}(h\gamma)^\alpha\Bigg],
\end{array}
$$
where $[\frac{m}{2}]$ is the integer part of $\frac{m}{2}$.

Thus when $\beta\leq 0$ we get
\begin{equation}
v(h\beta)=-D_1\sin(h\omega\beta)-D_2\cos(h\omega\beta)-Q_{m - 3}(h\beta), \label{eq.(3.8)}
\end{equation}
where
\begin{equation}
D_1=\frac{(-1)^m}{4\omega^{2m-1}}\sum\limits_{\gamma=0}^NC_{\gamma}(h\omega\gamma)\sin(h\omega\gamma),\ \ \ D_2=\frac{(-1)^m}{4\omega^{2m-1}}\sum\limits_{\gamma=0}^NC_{\gamma}(h\omega\gamma)\cos(h\omega\gamma),\label{eq.(3.9)}
\end{equation}
and
\begin{equation}
Q_{m - 3}(h\beta
)=\frac{(-1)^m}{2\omega^{2m-1}}\sum\limits_{k=[\frac{m}{2}]}^{m-2}
\sum\limits_{\alpha=m-2}^{2k-1}\frac{(-1)^{k+\alpha}(m-k-1)\omega^{2k-1}}{(2k-1-\alpha)!\ \alpha!}
(h\beta)^{2k-1-\alpha} \sum\limits_{\gamma=0}^NC_{\gamma}(h\gamma)^\alpha\label{eq.(3.10)}
\end{equation}
is a unknown polynomial of degree $m-3$ of $(h\beta)$.

Similarly, in the case $\beta\geq N$ for the convolution
$v(h\beta) = G_m(h\beta)*C_\beta$ we obtain
\begin{equation}
v(h\beta )= D_1\sin(h\omega\beta)+D_2\cos(h\omega\beta)+Q_{m - 3}(h\beta), \label{eq.(3.11)}
\end{equation}

We denote
\begin{eqnarray}
R_{m - 3}^{-}(h\beta )= R_{m - 3}(h\beta ) - Q_{m - 3} (h\beta), \ \ d_1^-=d_1-D_1,\ \ d_2^-=d_2-D_2, \label{eq.(3.12)}\\
R_{m - 3}^{+}(h\beta ) = R_{m - 3} (h\beta ) + Q_{m - 3}(h\beta),\ \ d_1^+=d_1+D_1,\ \ d_2^+=d_2+D_2,\label{eq.(3.13)}
\end{eqnarray}
where
$R_{m-3}^-(h\beta)=\sum_{\alpha=0}^{m-3}r_\alpha^-\cdot(h\beta)^\alpha,$
$R_{m-3}^+(h\beta)=\sum_{\alpha=0}^{m-3}r_\alpha^+\cdot(h\beta)^\alpha$.

Taking into account (\ref{eq.(3.6)}), (\ref{eq.(3.8)}) and
(\ref{eq.(3.11)}) we get the following problem

\begin{problem}\label{Pr3}
Find the solution of the equation
\begin{equation}
D_m(h\beta )*u(h\beta ) = 0,\,\,\,\,\,\,\,\,h\beta  \notin
[0,1]\label{eq.(3.14)}
\end{equation}
having the form:
\begin{equation}
u(h\beta ) = \left\{
\begin{array}{ll}
d_1^-\sin(h\omega\beta)+d_2^-\cos(h\omega\beta)+R_{m-3}^{-}(h\beta),& \beta\leq 0, \\[2mm]
\varphi(h\beta ),& 0 \le \beta  \le N, \\[2mm]
d_1^+\sin(h\omega\beta)+d_2^+\cos(h\omega\beta)+R_{m-3}^{+}(h\beta),&\beta\geq N.\\
\end{array}
\right.\label{eq.(3.15)}
\end{equation}
Here $R_{m-3}^{-}(h\beta)$ and $R_{m-3}^{+}(h\beta)$ are unknown
polynomials of degree $m-3$ with respect to $h\beta$.
\end{problem}

If we find $d_1^-,\ d_1^+,\ d_2^-,\ d_2^+$ and polynomials $R_{m-3}^{-}(h\beta)$, $R_{m-3}^{+}(h\beta)$ then from
(\ref{eq.(3.12)}), (\ref{eq.(3.13)}) we have
\begin{equation}
\label{eq.(3.16)}
\begin{array}{lll}
R_{m-3}(h\beta)=\frac{1}{2}\left(R_{m-3}^{+}(h\beta)+R_{m-3}^{-}(h\beta)\right),& d_k=\frac{1}{2}(d_k^++d_k^-),& k=1,2,\\
Q_{m-3}(h\beta)=\frac{1}{2}\left(R_{m-3}^{+}(h\beta)-R_{m-3}^{-}(h\beta)\right),& D_k=\frac{1}{2}(d_k^+-d_k^-),& k=1,2.
\end{array}
\end{equation}

Unknowns $d_1^-,\ d_1^+,\ d_2^-,\ d_2^+$ and polynomials $R_{m-3}^{-}(h\beta)$, $R_{m-3}^{+}(h\beta)$  can be found
from equation (\ref{eq.(3.14)}), using the function $D_m(h\beta)$
defined by (\ref{eq.(2.7)}). Then we obtain explicit form of the
function $u(h\beta )$ and from (\ref{eq.(3.7)}) we find the
coefficients $C_\beta$. Furthermore from (\ref{eq.(3.16)})  we get
$R_{m-3}(h\beta)$, $d_1$ and $d_2$.

Thus Problem \ref{Pr3} and respectively Problems \ref{Pr2} and
\ref{Pr1} will be solved.

In the next section we apply this algorithm to compute the
coefficients $C_{\beta}$, $\beta=0,1,...,N$, $d_1$, $d_2$ and $r_\alpha$,
$\alpha=0,1,...,m-3$ of the interpolation spline (\ref{eq.(1.4)})
for any $m\geq 2$ and $N+1\geq m$.

\section{Computation of coefficients of interpolation spline (\ref{eq.(1.4)})}
\setcounter{equation}{0} \setcounter{theorem}{0}
\setcounter{definition}{0} \setcounter{lemma}{0}
\setcounter{example}{0} \setcounter{figure}{0}
\setcounter{table}{0}

In this section, using the above algorithm, we obtain the explicit
formulas for the coefficients of the interpolation spline
(\ref{eq.(1.4)}) which, as we have proved in the previous section,
is the solution of Problem \ref{Pr1}.

It should be noted that the interpolation spline (\ref{eq.(1.4)}), the solution of Problem \ref{Pr1},
is exact for any polynomials of degree $m-3$ and for trigonometric functions $\sin\omega x$ and $\cos\omega x$.

In the sequel, we obtain the exact formulas for the coefficients of
the interpolation spline (\ref{eq.(1.4)}).  The result is the
following

\begin{theorem}\label{THM:4.1}
Coefficients of the interpolation spline (\ref{eq.(1.4)}), with
equally spaced nodes in the space $K_2(P_m)$, have the
following form
\begin{eqnarray*}
C_0&=&p\left[C\varphi(0)+\varphi(h)-d_1^-\sin(h\omega)+d_2^-\cos(h\omega)+\sum\limits_{\alpha=0}^{m-3}r_\alpha^-\cdot(-h)^\alpha\right]\\
  &&\qquad +\sum\limits_{k=1}^{m-1}\frac{A_kp}{\lambda_k}
\left[\sum\limits_{\gamma=0}^N\lambda_k^{\gamma}\varphi(h\gamma)+
M_k+\lambda_k^NN_k\right],\\
C_{\beta}&=&p\bigg[\varphi(h\beta-h)+C\varphi(h\beta)+\varphi(h\beta+h)\bigg]\\
&&\qquad +\sum\limits_{k=1}^{m-1}\frac{A_kp}{\lambda_k}
\left[\sum\limits_{\gamma=0}^N\lambda_k^{|\beta-\gamma|}\varphi(h\gamma)
+\lambda_k^{\beta}M_k+\lambda_k^{N-\beta}N_k\right],\\
&& \beta=1,2,...,N-1,\\
C_N&=&p\left[C\varphi(1)+\varphi(1-h)+d_1^+\sin(\omega+h\omega)+d_2^+\cos(\omega+h\omega)+\sum\limits_{\alpha=0}^{m-3}r_\alpha^+\cdot(1+h)^\alpha\right]\\
&&+\sum\limits_{k=1}^{m-1}\frac{A_kp}{\lambda_k}\left[\sum\limits_{\gamma=0}^N\lambda_k^{N-\gamma}\varphi(h\gamma)
+\lambda_k^NM_k+N_k\right],\\
d_k&=&\frac{1}{2}\left(d_k^++d_k^-\right),\ \ k=1,2,\\
r_\alpha&=&\frac{1}{2}\left(r_\alpha^++r_\alpha^-\right),\ \ \alpha=0,1,...,m-3,
\end{eqnarray*}
where
\begin{eqnarray}
M_k&=&\frac{\lambda_k[d_2^-(\cos(h\omega)-\lambda_1)-d_1^-\sin(h\omega)]}
{\lambda_k^2+1-2\lambda_k\cos(h\omega)}\nonumber\\
&&+
\sum\limits_{\alpha=1}^{m-3}r_\alpha^-(-h)^\alpha\sum\limits_{i=1}^\alpha\frac{\lambda_k^i\Delta^i
0^\alpha}{(1-\lambda_k)^{i+1}}+\frac{r_0^-\lambda_k}{1-\lambda_k},\label{eq.(4.1)}\\
N_k&=&\frac{\lambda_k[d_2^+(\cos(\omega+h\omega)-\lambda_k\cos
\omega)+d_1^+(\sin(\omega+h\omega)-\lambda_k\sin\omega)]}
{\lambda_k^2+1-2\lambda_k\cos(h\omega)}\nonumber\\
&&+\sum\limits_{\alpha=1}^{m-3}r_\alpha^+\left(\sum\limits_{j=1}^\alpha
C_\alpha^j h^j\sum\limits_{i=1}^j\frac{\lambda_k^i\Delta^i
0^j}{(1-\lambda_k)^{i+1}}+\frac{\lambda_k}{1-\lambda_k}\right)+\frac{r_0^+\lambda_k}{1-\lambda_k}\label{eq.(4.2)}
\end{eqnarray}
and $p$, $C$, $A_k$ are  defined by (\ref{eq.(2.8)}),(\ref{eq.(2.9)}), $\lambda_k$ are
the roots of the polynomial (\ref{eq.(2.11)}), $|\lambda_k|<1$,
$\Delta^i 0^\alpha=\sum_{l=1}^i (-1)^{i-l}C_i^ll^\alpha$, and $d_k^-$, $d_k^+$, $k=1,2,$
$r_\alpha^-,\ r_\alpha^+,$ $\alpha=0,1,...,m-3$, are defined from
the system (\ref{eq.(4.3)}), (\ref{eq.(4.4)}), (\ref{eq.(4.6)}),
(\ref{eq.(4.7)}).
\end{theorem}

\emph{Proof.} First we find the expressions for $d_2^-$ and
$d_2^+$. When $\beta=0$ and $\beta=N$ from (\ref{eq.(3.15)}) for
$d_2^-$ and $d_2^+$ we get
\begin{eqnarray}
d_2^-&=& \varphi(0)-r_0^-,\label{eq.(4.3)}\\
d_2^+&=&\frac{\varphi(1)}{\cos\omega}-d_1^+\tan\omega-\frac{1}{\cos\omega}\sum\limits_{\alpha=0}^{m-3}r_\alpha^+.\label{eq.(4.4)}
\end{eqnarray}
Now we have $2m-2$ unknowns $d_1^-$, $d_1^+$, $r_\alpha^-$, $r_\alpha^+$,
$\alpha=0,1,...,m-3$.

From equation (\ref{eq.(3.10)}), by choosing
$\beta=-1,-2,...,-(m-1)$ and $\beta=N+1,N+2,...,N+m-1$, we are
able to solve the previous system.

Taking into account  (\ref{eq.(3.15)}),
(\ref{eq.(4.3)}) and (\ref{eq.(4.4)}), from (\ref{eq.(3.14)}) we
get the following system
\begin{eqnarray}
&&-d_1^-\left[\sum_{\gamma=1}^{\infty}D_m(h\beta+h\gamma)\sin(h\omega\gamma)\right]+\sum_{\alpha=1}^{m-3}r_{\alpha}^-\left[(-h)^{\alpha}\sum_{\gamma=1}^{\infty}D_m(h\beta+h\gamma)
\gamma^{\alpha}\right]\nonumber \\
&&\quad +r_0^-\left[\sum_{\gamma=1}^{\infty}D_m(h\beta+h\gamma)(1-\cos(h\omega\gamma))\right]+d_1^+\left[
\sum_{\gamma=1}^{\infty}D_m(h(N+\gamma)-h\beta)\frac{\sin(h\omega \gamma)}
{\cos\omega}\right]\nonumber\\
&&\quad +\sum_{\alpha=1}^{m-3}r_{\alpha}^+\left[\sum_{j=1}^{\alpha}C_{\alpha}^jh^j\sum_{\gamma=1}^{\infty}D_m(h(N+\gamma)-h\beta)\gamma^j
+\sum_{\gamma=1}^{\infty}D_m(h(N+\gamma)-h\beta)\frac{\cos\omega-\cos(\omega+h\omega\gamma)}{\cos\omega}\right]\nonumber \\
&&\quad +r_0^+\left[\sum_{\gamma=1}^{\infty}D_m(h(N+\gamma)-h\beta)\frac{\cos\omega-\cos(\omega+h\omega\gamma)}{\cos\omega}\right]\nonumber \\
&&\qquad\qquad =-\sum_{\gamma=0}^ND_m(h\beta-h\gamma)\varphi(h\gamma)-\varphi(0)\left[\sum_{\gamma=1}^{\infty}D_m(h\beta+h\gamma)\cos(h\omega\gamma)\right]\nonumber\\
&&\qquad\qquad \qquad -\frac{\varphi(1)}{\cos\omega}\left[\sum_{\gamma=1}^{\infty}D_m(h(N+\gamma)-h\beta)\cos(\omega+h\omega\gamma)\right],\label{eq.(4.5)}
\end{eqnarray}
where $\beta=-1,-2,...,-(m-1)$ and $\beta=N+1,N+2,...,N+m-1$.

Now we consider the cases $\beta=-1,-2,...,-(m-1)$. From
(\ref{eq.(4.5)}) replacing $\beta$ by $-\beta$ and using
(\ref{eq.(2.7)}) and (\ref{eq.(2.5)}), after some calculations for
$\beta=1,2,...,m-1$, we get the following system of $m-1$ linear
equations
\begin{equation}
d_1^-B_{\beta}^-+\sum\limits_{\alpha=0}^{m-3}r_\alpha^-B_{\beta
\alpha}^-+d_1^+B_{\beta}^++\sum\limits_{\alpha=0}^{m-3}r_\alpha^+B_{\beta
\alpha}^+=T_\beta, \ \ \beta=1,2,...,m-1,\label{eq.(4.6)}
\end{equation}
where
\begin{eqnarray*}
B_{\beta}^-&=&-\left[\sum_{k=1}^{m-1}\frac{A_k}{\lambda_k}\sum_{\gamma=1}^{\infty}\lambda_k^{|\beta-\gamma|}\sin(h\omega\gamma)
+\sin(h\omega(\beta-1))+C\sin(h\omega\beta)+\sin(h\omega(\beta+1))\right],\\
B_{\beta\alpha}^-&=&(-h)^\alpha\left[\sum_{k=1}^{m-1}\frac{A_k}{\lambda_k}\sum_{\gamma=1}^{\infty}\lambda_k^{|\beta-\gamma|}\gamma^{\alpha}+
(\beta-1)^{\alpha}+C\beta^{\alpha}+(\beta+1)^{\alpha}\right],\\
B_{\beta0}^-&=&\sum_{k=1}^{m-1}\frac{A_k}{\lambda_k}\sum_{\gamma=1}^{\infty}\lambda_k^{|\beta-\gamma|}(1-\cos(h\omega\gamma))+
(1-\cos(h\omega(\beta-1)))\\
&&\qquad\qquad\qquad+C(1-\cos(h\omega\beta))+(1-\cos(h\omega(\beta+1))),\\
B_{\beta}^+&=&\frac{1}{\cos\omega}\sum_{k=1}^{m-1}\frac{A_k\lambda_k^{N+\beta}\sin(h\omega)}{\lambda_k^2+1-2\lambda_k\cos(h\omega)},\\
\end{eqnarray*}
\begin{eqnarray*}
B_{\beta\alpha}^+&=&\sum_{k=1}^{m-1}\frac{A_k\lambda_k^{N+\beta}}{\lambda_k}\Bigg[\sum_{j=1}^{\alpha}C_\alpha^jh^j\sum_{i=1}^j
\frac{\lambda_k^i\Delta^i0^j}{(1-\lambda_k)^{i+1}}+\frac{\lambda_k}{1-\lambda_k}\\
&&\qquad\qquad -\frac{\lambda_k[\cos(\omega+h\omega)-\lambda_k\cos(h\omega)]}{\cos \omega\ [\lambda_k^2+1-2\lambda_k\cos(h\omega)]}\Bigg],\\
B_{\beta0}^+&=&\sum_{k=1}^{m-1}A_k\lambda_k^{N+\beta}\Bigg[\frac{1}{1-\lambda_k}-\frac{\cos(\omega+h\omega)-\lambda_k\cos(h\omega)}
{\cos \omega\ [\lambda_k^2+1-2\lambda_k\cos(h\omega)]}\Bigg],\\
T_\beta&=&-\sum_{k=1}^{m-1}A_k\lambda_k^{\beta-1}\sum_{\gamma=0}^N\lambda_k^\gamma\varphi(h\gamma)-\varphi(0)\Bigg[\sum_{k=1}^{m-1}\frac{A_k}{\lambda_k}
\sum_{\gamma=1}^{\infty}\lambda_k^{|\beta-\gamma|}\cos(h\omega\gamma)\\
&&+\cos(h\omega(\beta-1))+C\cos(h\omega\beta)+\cos(h\omega(\beta+1))\Bigg]\\
&&\qquad\qquad-\frac{\varphi(1)}{\cos\omega}
\sum_{k=1}^{m-1}\frac{A_k\lambda_k^{N+\beta}[\cos(\omega+h\omega)-\lambda_k\cos\omega]}{\lambda_k^2+1-2\lambda_k\cos(h\omega)}.
\end{eqnarray*}
Here $\beta=1,2,...,m-1$ and $\alpha=1,2,...,m-3$.

Further, in (\ref{eq.(4.5)}), we consider the cases $\beta=N+1,N+2,...,N+m-1$. From
(\ref{eq.(4.5)}) replacing $\beta$ by $N+\beta$ and using
(\ref{eq.(2.7)}) and (\ref{eq.(2.5)}),  after some calculations
for $\beta=1,2,...,m-1$ we get the following system of $m-1$
linear equations
\begin{equation}
d_1^-A_{\beta}^-+\sum\limits_{\alpha=0}^{m-3}r_\alpha^-A_{\beta
\alpha}^-+d_1^+A_{\beta}^++\sum\limits_{\alpha=0}^{m-3}r_\alpha^+A_{\beta
\alpha}^+=S_\beta, \ \ \beta=1,2,...,m-1,\label{eq.(4.7)}
\end{equation}
where
\begin{eqnarray*}
A_{\beta}^-&=&-\sum_{k=1}^{m-1}\frac{A_k\lambda_k^{N+\beta}\sin(h\omega)}{\lambda_k^2+1-2\lambda_k\cos(h\omega)},\\
A_{\beta\alpha}^-&=&(-h)^\alpha\sum_{k=1}^{m-1}A_k\lambda_k^{N+\beta-1}\sum_{i=1}^{\alpha}\frac{\lambda_k^i\Delta^i0^\alpha}{(1-\lambda_k)^{i+1}},\\
A_{\beta0}^-&=&\sum_{k=1}^{m-1}\frac{A_k\lambda_k^{N+\beta}(\lambda_k+1)(\cos(h\omega)-1)}{(\lambda_k-1)(\lambda_k^2+1-2\lambda_k\cos(h\omega))}, \\
A_{\beta}^+&=&\frac{1}{\cos\omega}\Bigg[\sum_{k=1}^{m-1}\frac{A_k}{\lambda_k}\sum_{\gamma=1}\lambda_k^{|\beta-\gamma|}\sin(h\omega\gamma)\\
   &&+\sin(h\omega(\beta-1))+C\sin(h\omega\beta)+\sin(h\omega(\beta+1))\Bigg],\\
A_{\beta\alpha}^+&=&\sum_{j=1}^{\alpha}C_{\alpha}^jh^j\left[\sum_{k=1}^{m-1}\frac{A_k}{\lambda_k}\sum_{\gamma=1}^{\infty}\lambda_k^{|\beta-\gamma|}
\gamma^j+(\beta-1)^j+C\beta^j+(\beta+1)^j\right]\\
&&+\sum_{k=1}^{m-1}\frac{A_k}{\lambda_k}\sum_{\gamma=1}^{\infty}\lambda_k^{|\beta-\gamma|}+2+C-\frac{1}{\cos\omega}\Bigg[
\sum_{k=1}^{m-1}\frac{A_k}{\lambda_k}\sum_{\gamma=1}^{\infty}\lambda_k^{|\beta-\gamma|}\cos(\omega+h\omega\gamma)\\
\end{eqnarray*}
\begin{eqnarray*}
&&+\cos(\omega+h\omega(\beta-1))+C\cos(\omega+h\omega\beta)+\cos(\omega+h\omega(\beta+1))\Bigg],\\
A_{\beta0}^+&=&\sum_{k=1}^{m-1}\frac{A_k}{\lambda_k}\sum_{\gamma=1}^{\infty}\lambda_k^{|\beta-\gamma|}+2+C-\frac{1}{\cos\omega}\Bigg[
\sum_{k=1}^{m-1}\frac{A_k}{\lambda_k}\sum_{\gamma=1}^{\infty}\lambda_k^{|\beta-\gamma|}\cos(\omega+h\omega\gamma)\\
&&+\cos(\omega+h\omega(\beta-1))+C\cos(\omega+h\omega\beta)+\cos(\omega+h\omega(\beta+1))\Bigg],\\
S_\beta&=&-\sum_{k=1}^{m-1}\frac{A_k}{\lambda_k}\sum_{\gamma=0}^N\lambda_k^{N+\beta-\gamma}\varphi(h\gamma)
-\varphi(0)\sum_{k=1}^{m-1}\frac{A_k\lambda_k^{N+\beta}(\cos(h\omega)-\lambda_k)}{\lambda_k^2+1-2\lambda_k\cos(h\omega)}\\
&&-\frac{\varphi(1)}{\cos\omega}\Bigg[\sum_{k=1}^{m-1}\frac{A_k}{\lambda_k}\sum_{\gamma=1}^{\infty}
\lambda_k^{|\beta-\gamma|}\cos(\omega+h\omega\gamma)\\
&&+\cos(\omega+h\omega(\beta-1))+C\cos(\omega+h\omega\beta)+\cos(\omega+h\omega(\beta+1))\Bigg].
\end{eqnarray*}
Here $\beta=1,2,...,m-1$ and $\alpha=1,2,...,m-3$.

Thus for the unknowns $d_1^-$, $d_1^+$, $r_\alpha^-$, $r_\alpha^+$, $\alpha=0,1,...,m-3$ we have obtained system
(\ref{eq.(4.6)}), (\ref{eq.(4.7)}) of $2m-2$ linear equations.
Since our interpolation problem has a unique solution, the main
matrix of this system is non singular. Unknowns $d_1^-$, $d_1^+$, $r_\alpha^-$, $r_\alpha^+$, $\alpha=0,1,...,m-3$ can be found from system
(\ref{eq.(4.6)}), (\ref{eq.(4.7)}). Then taking into account
(\ref{eq.(3.16)}), using (\ref{eq.(4.3)}) and (\ref{eq.(4.4)}) we
have
\begin{eqnarray*}
d_k&=&\frac{1}{2}\left(d_k^++d_k^-\right),\ \ k=1,2,\\
r_\alpha&=&\frac{1}{2}\left(r_\alpha^++r_\alpha^-\right),\ \ \alpha=0,1,...,m-3.
\end{eqnarray*}
Now we find the coefficients $C_\beta$, $\beta=0,1,...,N$.

From (\ref{eq.(3.6)}), taking into account (\ref{eq.(3.15)}), we
deduce
\begin{eqnarray*}
C_\beta&=&\sum\limits_{\gamma=0}^ND_m(h\beta-h\gamma)\varphi(h\gamma)\\
&&+\sum\limits_{\gamma=1}^{\infty}D_m(h\beta+h\gamma)
\left[d_1^-\sin(-h\omega\gamma)+d_2^-\cos(h\omega\gamma)+\sum_{\alpha=0}^{m-3}r_{\alpha}^-(-h\gamma)^{\alpha}\right]\\
&&+\sum\limits_{\gamma=1}^{\infty}D_m(h(N+\gamma)-h\beta)
\left[d_1^+\sin(\omega+h\omega\gamma)+d_2^+\cos(\omega+h\omega\gamma)+\sum_{\alpha=0}^{m-3}r_{\alpha}^+(1+h\gamma)^{\alpha}\right],
\end{eqnarray*}
where $\beta=0,1,...,N.$\\
From here, using (\ref{eq.(2.7)}) and formula (\ref{eq.(2.5)}),
taking into account (\ref{eq.(4.1)}) and (\ref{eq.(4.2)}),  after
some calculations we arrive at the expressions of the coefficients
$C_\beta$, $\beta=0,1,...,N$ which are given in the assertion of
the theorem.

Theorem \ref{THM:4.1} is proved. \qed

\begin{remark}
From Theorem \ref{THM:4.1}, when $m=2$, we get Theorem 7 of \cite{Hay14} and Theorem 3.1 of \cite{HayMilShad14},
and when $m=2$, $\omega=1$ we get Theorem 3.1 of the work \cite{HayMilShad13}.
\end{remark}

\section*{Acknowledgements}  The part of this work has been done in the University of Santiago de Compostela, Spain. The author thanks
the program Erasmus Mundus Action 2, Marco XXI for financial
support (project number: Lot10-20112572).

\end{document}